\documentclass[10pt,twoside]{article}

\usepackage[centertags]{amsmath}
\usepackage{amsfonts}
\usepackage{amssymb}
\usepackage{amsthm}
\usepackage{epsfig}
\usepackage{color}
\allowdisplaybreaks[4]

\definecolor{c20}{rgb}{0.,0.7,0.}
\definecolor{c30}{rgb}{0.,0.,1.}
\definecolor{c40}{rgb}{1,0.1,0.7}
\definecolor{c50}{rgb}{1,0,0}
\definecolor{c60}{rgb}{1,0.9,0.1}

\def\wE#1{\textcolor{c30}{#1}}

\def\xE#1{\textcolor{c30}{#1}}

\def\xE#1{#1}

\def\wE#1{#1}

\newcommand{\COM}[1]{}

\newcommand{\ABs}[1]{\biggl \lvert #1 \biggr \rvert}

\newcommand{\BQN}{\begin{eqnarray}}
\newcommand{\EQN}{\end{eqnarray}}
\newcommand{\BQNY}{\begin{eqnarray*}}
\newcommand{\EQNY}{\end{eqnarray*}}
\newtheorem{theo}{Theorem}[section]
\newtheorem{pro}[theo]{Proposition}

\newtheorem{lem}[theo]{Lemma}

\newtheorem{exa}[theo]{Example}

\newcommand{\BL}{\begin{lem}}
\newcommand{\EL}{\end{lem}}

\newcommand{\BT}{\begin{theo}}
\newcommand{\ET}{\end{theo}}
\newcommand{\BEX}{\begin{exa}}
\newcommand{\EEX}{\end{exa}}
\newcommand{\BP}{\begin{pro}}
\newcommand{\EP}{\end{pro}}

\def\IF{\infty}

\begin{document}
\def\ZHT{Z(s,t)}
\def\ZHTY{Z^*(s,t)}
\def\MBT{\xE{M_{1/2}(T)}}
\def\MHT{\xE{M_{H}(T)}}
\def\MHTY{\xE{M^*_{H}(T)}}
\def\MHTD{\xE{M_{H}^{(\delta)}(T)}}
\def\MBTU{\xE{M_{H}(T_u)}}
\def\MHTU{\xE{M_{H}(T_u)}}

\title{\bf \Large Limit laws on extremes of non-homogeneous Gaussian random fields\thanks{Research supported by
National Science Foundation of China
(No. 11501250), China Postdoctoral Science Foundation (No. 2016M600460) and Natural Science Foundation of Zhejiang Province of China (No. LQ14A010012
).}}
\author{ Zhongquan Tan\footnote{E-mail address:  tzq728@163.com }
\\
{\small\it  School of Mathematical Sciences, Zhejiang University, Hangzhou 310027, China}\\
{\small\it  College of Mathematics, Physics and Information
Engineering, Jiaxing University, Jiaxing 314001, PR China;}\\
}

\bigskip

\date{\today}
 \maketitle

{\bf Abstract:} In this paper, by using the tail asymptotics derived by D\c{e}bicki, Hashorva and Ji (Ann. Probab. 2016),
we prove the Gumbel limit laws for the maximum of a class of non-homogeneous Gaussian random fields.
As an application of the main results, we derive the Gumbel limit law for Shepp statistics of fractional Brownian motion
and Gaussian integrated process.



\section{Introduction}

The studies on the Gumbel limit law for Gaussian processes have a long history and can date back to Pickands (1969).
Suppose that $\{X(t): t\in[0,\infty)\}$ is a
stationary Gaussian process with the covariance function $r(t)$  satisfying the following condition:
\begin{eqnarray}
\label{eq1.1}
r(t) = 1 -|t|^{\alpha} + o(|t|^{\alpha}),\ \ t \rightarrow 0,\ \ \mbox{and}\ \ \ r(t)<1,\ \ t> 0
\end{eqnarray}
with $\alpha\in (0, 2]$.
It is well-known (see e.g.
Pickands (1969), Leadbetter et al. (1983)) that
if further the so-called Berman's condition holds as follows
\begin{eqnarray*}
r(t)\ln t\rightarrow 0,\ \ \mbox{as}\ \  t\rightarrow \infty
\end{eqnarray*}
then the Gumbel limit law
\begin{eqnarray}
\label{eqA}
P\left(a_{T}\left(\sup_{0\leq t\leq T}X(t)-b_{T}\right)\leq x\right)\rightarrow \exp(-e^{-x})
\end{eqnarray}
holds for any $x\in R$, as $T\rightarrow\infty$, where
$$a_{T}=\sqrt{2\ln T}, \ \ b_{T}=\sqrt{2\ln T}+\frac{\ln[(2\pi)^{-1/2}\mathcal{H}_{\alpha}(2\ln T)^{-1/2+1/\alpha}]}{\sqrt{2\ln T}}.
$$
Here $\mathcal{H}_{\alpha}$ denotes the Pickands constant given by
$$\mathcal{H}_{\alpha}=\lim_{\lambda\rightarrow\infty}\mathcal{H}_{\alpha}[0,\lambda]/\lambda\in (0,\infty)$$
with
$$\mathcal{H}_{\alpha}[0,\lambda]=E\exp\left(\max_{t\in[0,\lambda]}\sqrt{2}B_{\alpha}(t)-t^{\alpha}\right)$$
and $B_{\alpha}$  a fractional Brownian motion (fBm) with Hurst parameter $\alpha/2
\in (0,1]$, that is, a  zero mean Gaussian process with stationary
increments such that $EB_{\alpha}^{2}(t)=|t|^{\alpha}$. To derive the Gumbel limit law (\ref{eqA}), the following well-known Pickands asymptotics (see e.g.
Pickands (1969),  Berman (1974), Leadbetter et al. (1983)) plays a crucial role, i.e.,
\begin{eqnarray}
\label{eq1.3}
P\left(\sup_{t\in[0,T]}X(t)>u\right)=T\mathcal{H}_{\alpha}u^{2/\alpha}\Psi(u)(1+o(1)),
\end{eqnarray}
as $u\rightarrow\infty$ for some fixed $T\in(0,\infty)$,  where $\Psi(\cdot)$ denotes the  tail distribution of a standard normal
random variable.
For some recent work on the tail asymptotics for extremes,
we refer to Chan and Lai (2006), D\c{e}bicki, Hashorva, Ji and Tabi\'{s} (2015),  Cheng and Xiao (2016,2017) and the references therein.

The investigation of (\ref{eqA}) for Gaussian processes and general stochastic processes has received a lot of attention.
Mittal and Ylvisaker (1975) extended (\ref{eqA}) to the strongly dependent Gaussian case;
H\"{u}lser (1990) investigated (\ref{eqA}) for locally stationary Gaussian case, which is recently further extended to Gaussian random fields
on manifolds by Qiao and Polonic (2017). We refer to McCormick (1980), Konstant and  Piterbarg (1993) and
Piterbarg (1996) for further extensions to Gaussian processes and fields;
Leadbetter and Rootz\'{e}n (1982) and Albin (1990) for  stationary non-Gaussian processes.
For more related extensions, we refer to D\c{e}bicki, Hashorva, Ji and Ling (2015) and the reference therein.

In many applied fields, the Gumbel limit laws for extremes of Gaussian processes
 play a very important role. In approximation theory, Seleznjev (1991, 1996), H\"{u}lser (1999) and H\"{u}lser et al. (2003)
applied the Gumbel limit law for Gaussian processes to investigate  the deviation processes of
some piecewise linear interpolation problems;
In nonparametric statistics, the absolute deviations of many types of density estimators
obey the Gumbel limit law, see e.g. Bickel and Rosenblatt (1973) and Gin\'{e} et al. (2003).
In applied statistics, there are also many confidence intervals and bands,
which are constructed based on the Gumbel limit law of the estimators,
since extremes themselves are also type of very important estimators, see e.g. Gin\'{e} and Nickl (2010).
For some recent studies on applications of Gumbel limit laws, we refer to Sharpnack and Arias-Castro (2016) and Qiao and Polonik (2016).

Define
\begin{eqnarray}
\label{eq1.4}
X(s,t)=Y(s+t)-Y(s),
\end{eqnarray}
 where $\{Y(t),t\geq 0\}$ is a Gaussian process.
The process $\xi_{S}(t)=\sup_{0\leq s\leq S}X(s,t)$ is referred to as the Shepp
statistics in many recent works. Zholud (2008) studied the maximum of the process $\xi_{S}(t)$ and
established the Gumbel limit law when $Y(t)$ is a Brownian motion.
Hashorva and Tan (2013) and Tan and Yang (2014) extended the result to fractional Brownian motion.
We refer to Piterbarg (2001), H\"{u}lser and Piterbarg (2004a) and Hashorva et al. (2013) for related work on the fractional Brownian motion.

In this paper, we generalize model (\ref{eq1.4}) and impose directly some restrictions on the Gaussian random fields.
We first consider the Gumbel limit law for the process $\varsigma_{T}(s)=\sup_{0\leq t\leq T}X(s,t)$ for some fixed $T>0$,
where $X(s,t)$ is a type of non-homogeneous Gaussian random field. Then we use the obtained results to
derive the Gumbel limit law for Shepp statistics.
Noting that $\varsigma_{T}(t)$ is no longer Gaussian process, we can not derive the Gumbel limit laws from the Gaussian case directly.
However, $\varsigma_{T}(t)$ also doesn't satisfy the conditions imposed on general stochastic processes, such as those given by
Leadbetter and Rootz\'{e}n (1982) and Albin (1990).
We will follow the method used in Chapter 12 in Leadbetter et al. (1983).
The tail asymptotic result of extremes of the field $X(s,t)$  is a key tool, which has been derived by D\c{e}bicki et al. (2016).

The rest of the paper is organized as follows. In Section 2, we give some tail asymptotic results from D\c{e}bicki et al. (2016).
In Section 3, we state
the main results of the paper, and in Section 4, we present two applications. The technical
proofs are gathered in Section 5, while in Section 6 we give two auxiliary results.

\section{Preliminaries}
In this section, we present the tail asymptotic result provided by D\c{e}bicki et al. (2016).
Suppose that $\{X(s,t), (s,t)\in [0,\infty)\times[0,T]\}$ with fixed $T$ is a centered Gaussian random field with variance function and correlation function $\sigma^{2}(s,t)$ and
$r(s,t,s',t')$, respectively. Suppose the following assumptions hold.

{\bf Assumption A1:} there exists some positive function $\sigma(t)$ which attains its unique maximum on $[0,T]$ at fixed $T$, and further
$$\sigma(s,t)=\sigma(t),\ \  \forall (s,t)\in [0,\infty)\times[0,T],\ \ \sigma(t)=1-b(T-t)^{\beta}(1+o(1)),\ \ t\uparrow T$$
holds for some $\beta,b>0$.

{\bf Assumption A2:} there exist constants $a_{1}>0, a_{2}>0, a_{3}\neq 0$ and $\alpha_{1}, \alpha_{2}\in (0,2]$ such that
$$r(s,t,s',t')=1-(|a_{1}(s-s')|^{\alpha_{1}}+|a_{2}(t-t')+a_{3}(s-s')|^{\alpha_{2}})(1+o(1))$$
holds uniformly with respect to $s,s'\in [0,L]$ with some constant $L>0$ as $|s-s'|\rightarrow0$, $\min(t,t')\uparrow T$ and
further, there exists some constant $\delta_{0}\in (0,T)$ such that
$$r(s,t,s',t')<1$$
for any $s,s'\in [0,L]$ satisfying $s\neq s'$ and $t,t'\in[\delta_{0}, T]$.

{\bf Assumption A3:} There  exist positive constants $\gamma_{1}, \gamma_{2}, \gamma$ and $\mathcal{C}$ such that
$$E(X(s,t)-X(s',t'))^{2}\leq \mathcal{C} (|t-t'|^{\gamma}+|s-s'|^{\gamma})$$
holds for all $t,t'\in[\gamma_{1},T]$, $s,s'\in[0,L]$ satisfying $|s-s'|<\gamma_{2}$.

To state the tail asymptotics for the maximum of the field $X(s,t)$ under assumptions {\bf A1-A3}, we need the so-called Piterbarg constants and Pickands-Piterbarg constants, respectively.
The Piterbarg constant $\mathcal{P}_{\alpha}^{b}$ with constant $b>0$ is
defined as
$$\mathcal{P}_{\alpha}^{b}=\lim_{\lambda\rightarrow\infty}E\exp\left(\max_{t\in[0,\lambda]}\sqrt{2}B_{\alpha}(t)-(1+b)|t|^{\alpha}\right)\in (0,\infty).$$
For some constants $a_{1}>0, a_{2}>0, a_{3}\neq 0, b>0$, let
$$Y(s,t)=\widetilde{B}_{\alpha}(a_{1}s)+B_{\alpha}(a_{2}t-a_{3}s),\ \ \sigma_{Y}^{2}(s,t)=Var(Y(s,t))$$
and
$$\mathcal{H}_{Y}^{b}[\lambda_{1},\lambda_{2}]=E\exp\left(\max_{(s,t)\in[0,\lambda_{1}]\times[0,\lambda_{2}]}\sqrt{2}Y(s,t)-\sigma_{Y}^{2}(s,t)-b|t|^{\alpha}\right),$$
where $\widetilde{B}_{\alpha}$ and $B_{\alpha}$ are two independent fractional Brownian motions (fBms).
The Pickands-Piterbarg constant is defined as
$$\mathcal{M}_{Y,\alpha}^{b}=\lim_{\lambda_{1}\rightarrow\infty}\lim_{\lambda_{2}\rightarrow\infty}\frac{1}{\lambda_{1}}\mathcal{H}_{Y}^{b}[\lambda_{1},\lambda_{2}].$$

Under the above assumptions, D\c{e}bicki et al. (2016) derived the following result.

\BT\label{Th:main1} Let $\{X(s,t), (s,t)\in [0,\infty)\times[0,T]\}$ with fixed  $T$ be a centered Gaussian random field with a.s. continuous sample
paths. Suppose that assumptions {\bf A1-A3} are satisfied with the parameters mentioned therein, we have as $u\rightarrow\infty$,
$$P\left(\sup_{(s,t)\in [0,L]\times[0,T]}X(s,t)>u\right)=L\mu(u)(1+o(1)),$$
where for $\beta> \max\{\alpha_{1},\alpha_{2}\}$
$$\mu(u)=\Gamma(1/\beta+1)\prod_{k=1}^{2}(a_{k}\mathcal{H}_{\alpha_{k}})b^{-\frac{1}{\beta}}
u^{\frac{2}{\alpha_{1}}+\frac{2}{\alpha_{2}}-\frac{2}{\beta}}\Psi(u);$$
for $\beta=\alpha_{2}=\alpha_{1}$
$$\mu(u)=\mathcal{M}_{Y,\alpha_{1}}^{b}
u^{\frac{2}{\alpha_{1}}}\Psi(u);$$
for $\beta=\alpha_{2}>\alpha_{1}$
$$\mu(u)=a_{1}a_{2}\mathcal{P}_{\alpha_{2}}^{ba_{2}^{-\alpha_{2}}}\mathcal{H}_{\alpha_{1}}
u^{\frac{2}{\alpha_{1}}}\Psi(u);$$
for $\beta<\alpha_{2}=\alpha_{1}$
$$\mu(u)=(a_{1}^{\alpha_{1}}+|a_{3}|^{\alpha_{1}})^{\frac{1}{\alpha_{1}}}\mathcal{H}_{\alpha_{1}}
u^{\frac{2}{\alpha_{1}}}\Psi(u);$$
for $\beta<\alpha_{2}$ and $\alpha_{1}<\alpha_{2}$
$$\mu(u)=a_{1}\mathcal{H}_{\alpha_{1}}
u^{\frac{2}{\alpha_{1}}}\Psi(u);$$
for $\beta=\alpha_{1}>\alpha_{2}$
$$\mu(u)=a_{1}\mathcal{P}_{\alpha_{1}}^{b(\frac{|a_{3}|}{a_{1}a_{2}})^{\alpha_{1}}}\mathcal{H}_{\alpha_{2}}
u^{\frac{2}{\alpha_{2}}}\Psi(u);$$
for $\beta<\alpha_{1}$ and $\alpha_{2}<\alpha_{1}$
$$\mu(u)=|a_{3}|\mathcal{H}_{\alpha_{2}}
u^{\frac{2}{\alpha_{2}}}\Psi(u).$$
\ET
This result is very powerful since it can be used to derive the exact tail asymptotics for many type of statistics, such as
Shepp statistics for Gaussian processes, Brownian bridge and fBm, maximum loss and span of Gaussian processes, see D\c{e}bicki et al. (2016) for details.

\section{Main Result}

Note that the assumptions {\bf A1} and {\bf A2} are local conditions.
To derive the Gumbel limit law, we need to impose the following Berman-type weak dependence condition, which is a global condition.

{\bf Assumption A4:}  Assume that for $c=1+\varepsilon I(\beta\geq\max\{\alpha_{1},\alpha_{2}\})$ with some constant $\varepsilon>0$ the
function
$$\delta(v) := \sup\{|r(s,t,s',t')|, |s-s'|\ge v,s,s'\in [0,\infty), t,t'\in[0,T] \}$$ is such that
\BQN \label{BermD}
\lim_{v\to \IF} \delta(v) (\ln v)^c=0,
\EQN
where $I(\cdot)$ denotes the indicator function.

We state  now the main result.

\BT\label{Th:main23} Let $\{X(s,t), (s,t)\in [0,\infty)\times[0,T]\}$ with fixed $T$ be a centered Gaussian random field with a.s. continuous sample
paths. Suppose that assumptions {\bf A1-A4} are satisfied with the parameters mentioned therein.
In addition, assume that $\{X(s,t), (s,t)\in [0,\infty)\times[0,T]\}$ is homogeneous with respect to the first factor $s$.
Then
$$\lim_{S\to \infty}\sup_{x\in \mathbb{R}}
\ABs{P\left(a_{S} \left(\sup_{(s,t)\in [0,S]\times[0,T]}X(s,t)-b_{S}\right)\leq x\right)-\exp(-e^{-x})}=0,$$
where $a_{S}=\sqrt{2\ln S}$
$$\ \ b_{S}=a_{S}+a_{S}^{-1}\omega_{S}$$
with, for $\beta> \max\{\alpha_{1},\alpha_{2}\}$
$$\omega_{S}=\ln\left((2\pi)^{-1/2}\Gamma(1/\beta+1)\prod_{k=1}^{2}(a_{k}\mathcal{H}_{\alpha_{k}})b^{-\frac{1}{\beta}}
a_{S}^{\frac{2}{\alpha_{1}}+\frac{2}{\alpha_{2}}-\frac{2}{\beta}-1}\right);$$
for $\beta=\alpha_{2}=\alpha_{1}$
$$\omega_{S}=\ln\left((2\pi)^{-1/2}\mathcal{M}_{Y,\alpha_{1}}^{b}
a_{S}^{\frac{2}{\alpha_{1}}-1}\right);$$
for $\beta=\alpha_{2}>\alpha_{1}$
$$\omega_{S}=\ln\left((2\pi)^{-1/2}a_{1}a_{2}\mathcal{P}_{\alpha_{2}}^{ba_{2}^{-\alpha_{2}}}\mathcal{H}_{\alpha_{1}}
a_{S}^{\frac{2}{\alpha_{1}}-1}\right);$$
for  $\beta<\alpha_{2}=\alpha_{1}$
$$\omega_{S}=\ln\left((2\pi)^{-1/2}(a_{1}^{\alpha_{1}}+|a_{3}|^{\alpha_{1}})^{\frac{1}{\alpha_{1}}}\mathcal{H}_{\alpha_{1}}
a_{S}^{\frac{2}{\alpha_{1}}-1}\right);$$
for $\beta<\alpha_{2}$ and $\alpha_{1}<\alpha_{2}$
$$\omega_{S}=\ln\left((2\pi)^{-1/2}a_{1}\mathcal{H}_{\alpha_{1}}
a_{S}^{\frac{2}{\alpha_{1}}-1}\right);$$
for $\beta=\alpha_{1}>\alpha_{2}$
$$\omega_{S}=\ln\left((2\pi)^{-1/2}a_{1}\mathcal{P}_{\alpha_{1}}^{b(\frac{|a_{3}|}{a_{1}a_{2}})^{\alpha_{1}}}\mathcal{H}_{\alpha_{2}}
a_{S}^{\frac{2}{\alpha_{2}}-1}\right);$$
for $\beta<\alpha_{1}$ and $\alpha_{2}<\alpha_{1}$
$$\omega_{S}=\ln\left((2\pi)^{-1/2}|a_{3}|\mathcal{H}_{\alpha_{2}}
a_{S}^{\frac{2}{\alpha_{2}}-1}\right).$$
\ET

{\bf Remark 3.1:} Assumption A4 is a weakly dependent condition. If $\lim_{v\to \IF} \delta(v) (\ln v)^c=d>0$, then the field $X(s,t)$ will
possess some strongly dependent property with respect to the first parameter. In this case, the limit distribution will be no longer Gumbel distribution,
see Mittal and Ylvisaker (1975) and Tan et al. (2012) for some related results about strongly dependent Gaussian processes.

\section{Applications}

In this section, we give two applications of our main results.
We derive the exact tail asymptotics and Gumbel limit laws for Shepp statistics. The obtained results are of independent interest.

Throughout this section, let $\{X(t),t\geq0\}$ be a centered Gaussian process and define
$$Z(s,t)=X(s+t)-X(s),\ \ (s,t)\in[0,\infty)\times[0,T],$$
for some fixed $T>0$. The Shepp statistic $\sup_{0\leq s\leq S}Z(s,t)$ which was introduced by
Shepp (see Shepp 1966,1971) play a vary important role in statistics. Other important results for the
Shepp statistics can be found in
Cressie (1980), Deheuvels and Devroye (1987), Siegmund and Venkatraman (1995),  Dumbgen and Spokoiny (2001) and  Kabluchko (2011).
The limit properties of extremes of Shepp statistics when $X(t)$ is a fBm have been
studied by  Zholud (2008) and Hashorva and Tan (2013), Tan and Yang (2015) and Tan and Chen (2016).
Applying Theorem \ref{Th:main23}, we study the limit properties of extremes of Shepp statistics for a more general Gaussian process $X(t)$, which is a stationary Gaussian process or non-stationary Gaussian process with stationary increments.

\subsection{Stationary case}

Let $\{X(t),t\geq0\}$ be a centered stationary Gaussian process.
Suppose the covariance function $r_{X}$ of $\{X(t),t\geq0\}$ satisfies the following conditions:

{\bf Assumption B1:} $r_{X}(t)$ attains its minimum on $[0,T]$ at the unique point $T$;

{\bf Assumption B2:} there exist positive constants $\alpha_{1}, a_{1},a_{2}$ and $\alpha_{2}\in(0,2)$ such that
$$r_{X}(t)=r_{X}(T)+a_{1}|t-T|^{\alpha_{1}}(1+o(1)),\ \ t\rightarrow T,\ \ \mbox{and}\ \ \ r_{X}(t)=1-a_{2}t^{\alpha_{2}}(1+o(1)),\ \ t\rightarrow 0;$$
{\bf Assumption B3:} $r_{X}(s)<1$ for $s>0$.

For simplicity, write $\rho_{T}=\sqrt{2(1-r_{X}(T))}$ and $b_{i}=a_{i}/\rho_{T}^{2}$, $i=1,2$.
\BP\label{Pro:main3:1} Let $Z(s,t)$ be defined as above. Suppose that $r_{X}(t)$ satisfies conditions ${\bf B1-B3}$. In addition, suppose that
$r_{X}(t)$ is twice continuously differentiable on $[\tau,\infty)$ for some $\tau>0$ and the limit of twice derivative $\lim_{t\rightarrow T}|\ddot{r}_{X}(t)|\in(0,\infty)$.
Furthermore, if $\ddot{r}_{X}(t)(\ln t)^{c}=o(1)$ with $c=1+\varepsilon I(\alpha_{1}\geq\alpha_{2})$ and some constant $\varepsilon>0$ as $t\rightarrow\infty$, then
\BQN \label{thm3:1}
\lim_{S\to \infty}\sup_{x\in \mathbb{R}}
\ABs{P\left(a_{S}\left(\sup_{(s,t)\in[0,S]\times[0,T]}Z(s,t) -b_{S}\right)\leq x\right)-\exp\{-e^{-x}\}}=0,
\EQN
where $a_S=\rho_{T}\sqrt{2\ln S }$,
$$ b_{S}=a_{S}+a_{S}^{-1}\omega_{S}$$
with, for $\alpha_{1}>\alpha_{2}$
$$\omega_{S}=\ln (\Gamma(\frac{1}{\alpha_{1}}+1)\mathcal{H}_{\alpha_{2}}^{2}b_{2}^{\frac{2}{\alpha_{2}}}b_{1}^{-\frac{1}{\alpha_{1}}}
(2\pi)^{-1/2}a_{S}^{\frac{4}{\alpha_{2}}-\frac{2}{\alpha_{1}}-1});$$
for $\alpha_{1}=\alpha_{2}$
$$\omega_{S}=\ln (\mathcal{M}_{Y,\alpha_{1}}^{b_{1}}
(2\pi)^{-1/2}a_{S}^{\frac{2}{\alpha_{2}}-1})$$
with
$Y=Y(s,t)=\widetilde{B}_{\alpha_{2}}(b_{2}^{\frac{1}{\alpha_{2}}}s)+B_{\alpha_{2}}(b_{2}^{\frac{1}{\alpha_{2}}}t-b_{2}^{\frac{1}{\alpha_{2}}}s)$;
for $\alpha_{1}<\alpha_{2}$
$$\omega_{S}=\ln ((2b_{2})^{\frac{1}{\alpha_{2}}}\mathcal{H}_{\alpha_{2}}
(2\pi)^{-1/2}a_{S}^{\frac{2}{\alpha_{2}}-1}).$$
\EP

{\bf Example 4.1:} The Ornstein-Uhlenbeck process with covariance function $r_{X}(t)=e^{-|t|^{\alpha}}$
and the generalized Cauchy model with covariance function $r_{X}(t)=(1+|t|^{\alpha})^{-\beta}$ with $\alpha\in(0,2)$ and $\beta>0$ satisfy the conditions of Proposition \ref{Pro:main3:1}.

\subsection{Non-Stationary case}

Let $\{X(t),t\geq0\}$ be a centered non-stationary Gaussian process with stationary increment and variance function $\sigma_{X}^{2}(t)$, a.s. continuous sample paths. Recall that $X(t)$ is said to have stationary increments if the law of the process $\{X(t+t_{0})-X(t_{0}),t\in \mathbb{R}\}$ does not depend on
the choice of $t_{0}$.
To study the maximum of $Z(s,t)$, we only need to impose some conditions on the variogram $\gamma(t)=\mathbb{E}(X(t)-X(0))^{2}$ of $X$.
Note that for this case the variogram is $\gamma_{X}(t)=\sigma_{X}^{2}(t)$.
Suppose that the variance function $\sigma_{X}^{2}(t)$ of $\{X(t),t\geq0\}$ satisfies the following conditions:

{\bf Assumption C1:} $\sigma_{X}(t)$ attains its maximum on $[0,T]$ at the unique point $T$, and further
$$\sigma_{X}(t)=1-b(T-t)^{\beta}(1+o(1)),\ \ t\uparrow T$$
holds for some $\beta, b>0$.

{\bf Assumption C2:} $\sigma_{X}^{2}(t)$ is twice continuously differentiable on $[\tau,\infty)$ for $\tau>0$ with limit of twice derivative
 $\lim_{t\rightarrow T}|\ddot{\sigma}_{X}^{2}(t)|\in(0,\infty)$  and further
$$\sigma_{X}^{2}(t)=(at)^{\alpha}(1+o(1)),\ \ t\rightarrow 0$$
holds for some $\alpha\in(0,2], a>0$. 

{\bf Assumption C3:} $\ddot{\sigma}_{X}^{2}(t)(\ln t)^{c}\rightarrow 0$ with $c=1+\varepsilon I(\beta\geq\alpha)$ and some constant $\varepsilon>0$ as $t\rightarrow\infty$.

\BP\label{pro:main3:1} Let $Z(s,t)$ be defined as above. Suppose that $\sigma_{X}(t)$ satisfies conditions ${\bf C1,C2}$. We have for some constant $L>0$
 $$P\left(\sup_{(s,t)\in[0,L]\times[0,T]}Z(s,t)>u\right)=L\mu(u)(1+o(1)),$$
where for $\alpha<\beta$
$$\mu(u)=2^{-\frac{2}{\alpha}}\Gamma(1/\beta+1)a^{2}\mathcal{H}_{\alpha}^{2}b^{-\frac{1}{\beta}}
u^{\frac{4}{\alpha}-\frac{2}{\beta}}\Psi(u);$$
for $\alpha=\beta$
$$\mu(u)=\mathcal{M}_{Y}^{b}u^{\frac{2}{\alpha}}\Psi(u)$$
with
 $Y=Y(s,t)=\widetilde{B}_{\alpha}(2^{-1/\alpha}as)+B_{\alpha}(2^{-1/\alpha}at-2^{-1/\alpha}as)$;
for $\alpha>\beta$
$$\mu(u)=a\mathcal{H}_{\alpha}u^{\frac{2}{\alpha}}\Psi(u).$$
Furthermore, if condition ${\bf C3}$ holds, then
\BQN \label{thm3:1}
\lim_{S\to \infty}\sup_{x\in \mathbb{R}}
\ABs{P\left(a_{S}\left(\sup_{(s,t)\in[0,S]\times[0,T]}Z(s,t) -b_{S}\right)\leq x\right)-\exp\{-e^{-x}\}}=0,
\EQN
where $a_S=\sqrt{2\ln S }$, and
$$ b_{S}=a_{S}+a_{S}^{-1}\omega_{S}$$
with for $\alpha<\beta$
$$\omega_{S}=\ln (2^{-\frac{2}{\alpha}}\Gamma(1/\beta+1)a^{2}\mathcal{H}_{\alpha}^{2}b^{-\frac{1}{\beta}}
(2\pi)^{-1/2}a_{S}^{\frac{4}{\alpha}-\frac{2}{\beta}-1});$$
for $\alpha=\beta$
$$\omega_{S}=\ln (\mathcal{M}_{Y}^{b}
(2\pi)^{-1/2}a_{S}^{\frac{2}{\alpha}-1});$$
for $\alpha>\beta$
$$\omega_{S}=\ln (a\mathcal{H}_{\alpha}(2\pi)^{-1/2}a_{S}^{\frac{2}{\alpha}-1}).$$
\EP

We illustrate Proposition \ref{pro:main3:1} by the following two examples on the fBm and Gaussian integrated process.

{\bf Example 4.2:}  Let $B_{H_{i}}(t)$, $i=1,2,\ldots,n$ be a sequence of independent fBms with Hurst index $H_{i}\in(0,1)$ and $\lambda_{i}$ be a positive
sequence satisfying $\sum_{i=1}^{n}\lambda_{i}^{2}=1$.
Since given $H=H_{1}=H_{2}$ we have $\lambda_{1}B_{H_{1}}(t)+\lambda_{2}B_{H_{2}}(t)=^{d}\sqrt{\lambda_{1}^{2}+\lambda_{2}^{2}}B_{H}(t)$, we suppose that
$$H:=H_{1}<H_{2}<\cdots<H_{n}.$$
Let
$X(t)=\sum_{i=1}^{n}\lambda_{i}T^{-1/2}B_{H_{i}}(t)$
and $Z(s,t)$ be defined as above. We have for some constant $L>0$
$$P\left(\sup_{(s,t)\in[0,L\emph{}]\times[0,T]}Z(s,t) >u\right)=L\mu(u)(1+o(1)),$$
as $u\rightarrow\infty$, where
for $H\in(0,1/2)$
$$\mu(u)=2^{-\frac{1}{H}}\mathcal{H}_{2H}^{2}\lambda_{1}^{\frac{2}{H}}(\sum_{i=1}^{n}\lambda_{i}^{2}H_{i})^{-1}u^{\frac{2}{H}-2}\Psi(u);$$
for $H=1/2$
$$\mu(u)=\mathcal{M}_{Y,1}^{\frac{1}{2}}u^{2}\Psi(u)$$
with
$Y=Y(s,t)=\widetilde{B}_{1}(2^{-1}s)+B_{1}(2^{-1}(t-s))$;
for $H\in (1/2,1)$
$$\mu(u)=\lambda_{1}^{\frac{1}{H}}\mathcal{H}_{2H}u^{\frac{1}{H}}\Psi(u)$$
and
\BQN \label{thm3:2}
\lim_{S\to \infty}\sup_{x\in \mathbb{R}}
\ABs{P\left(a_{S}\left(\sup_{(s,t)\in[0,S]\times[0,T]}Z(s,t) -b_{S}\right)\leq x\right)-\exp\{-e^{-x}\}}=0,
\EQN
where $a_S=\sqrt{2 \ln S }$, and $$b_S=a_S+ a_S^{-1}\omega_{S}$$
with for $H\in (0,1/2)$
$$\omega_{S}= \ln\Bigl((2\pi)^{-1/2}2^{-\frac{1}{H}}\mathcal{H}_{2H}^{2}\lambda_{1}^{\frac{2}{H}}(\sum_{i=1}^{n}\lambda_{i}^{2}H_{i})^{-1}a_{S}^{\frac{2}{H}-3}\Bigr)$$
for any $H=1/2$
$$\omega_{S}=\ln\Bigl(2\pi^{-1/2})\mathcal{M}_{Y,1}^{\frac{1}{2}}a_{S}\Bigr)$$
and for $H\in (1/2,1)$
$$\omega_S=\ln\Bigl((2\pi)^{-1/2})\lambda_{1}^{\frac{1}{H}}\mathcal{H}_{2H}a_{S}^{\frac{1}{H}-1}\Bigr).$$

Next, we consider the Gaussian integrated process. For related studies, we refer to D\c{e}bicki (2002) and H\"{u}sler and Piterbarg (2004b).

{\bf Example 4.3:} Let $\{\zeta(t),t\geq0\}$ be a centered stationary Gaussian process with variance one and
suppose the covariance function $r_{\zeta}(t)$ of $\{\zeta(t),t\geq0\}$ satisfying the following conditions:

{\bf Assumption D1:} $r_{\zeta}(t)\in C([0,\infty))$ and $\int_{0}^{t}r_{\zeta}(s)ds>0$ for $t\in(0,T]$;

{\bf Assumption D2:} $r_{\zeta}(t)=1-t^{\theta}(1+o(1))$ as $t\rightarrow 0^{+}$ with $\theta\in(0,2]$;

{\bf Assumption D3:} $r_{\zeta}(t)\ln t=o(1)$ as $t\rightarrow\infty$.

Define Gaussian integrated processes as
$X(t)=\int_{0}^{t}\zeta(s)ds$
and let $Z(s,t)$ be defined as above.

 If conditions {\bf D1,D2} are satisfied, we have for some constant $L>0$
$$P\left(\sup_{(s,t)\in[0,L]\times[0,T]}Z(s,t) >u\right)=L\frac{1}{\sqrt{\pi}}
u\Psi(u)(1+o(1)),$$
as $u\rightarrow\infty$.
If further condition {\bf D3} holds, we have
\BQN \label{thm3:2}
\lim_{S\to \infty}\sup_{x\in \mathbb{R}}
\ABs{P\left(a_{S}\left(\sup_{(s,t)\in[0,S]\times[0,T]}Z(s,t) -b_{S}\right)\leq x\right)-\exp\{-e^{-x}\}}=0,
\EQN
where $a_S=\sqrt{2\ln S }$, and
$$ b_S=a_S+ a_S^{-1}  \ln\Bigl((2)^{-1/2}\pi^{-1} \Bigr).$$

\section{Proofs}

We need the following lemmas to prove Theorem \ref{Th:main23}. For simplicity, write $u=u_{S}(x)=a_{S}^{-1}x+b_{S}$ in the following part.

\BL\label{Le:A1} Let $\delta_{u}=u^{-2/\beta}(\ln u)^{2/\beta}$.
Under the conditions of Theorem 3.1,  we have for some constant $S_{0}>0$
$$\bigg|P\left(\sup_{(s,t)\in [0,S_{0}]\times[0,T]}X(s,t)\leq u\right)
-P\left(\sup_{(s,t)\in [0,S_{0}]\times[T-\delta_{u},T]}X(s,t)\leq u\right)\bigg|\bigg/P\left(\sup_{(s,t)\in [0,S_{0}]\times[0,T]}X(s,t)> u\right)\rightarrow 0$$
as $u\rightarrow\infty$.
\EL

\textbf{Proof:} It can be found in the proof of Theorem 2.2 of D\c{e}bicki et
al. (2016). $\Box$

For given $\epsilon> 0$, we divide interval $[0, S]$ into intervals of length 1, and split each of them onto subintervals $I_{j}^{\epsilon}$, $I_{j}$ of length $\epsilon$, $1-\epsilon$, $j=1,2,\cdots,\lfloor S\rfloor,$ respectively, where $\lfloor x\rfloor$ denotes the integral part of $x$. It can be easily seen that a possible remaining interval with length smaller than 1 plays no role in our consideration. We denote this interval by $J$.

\BL\label{Le:A2} Under the conditions of Theorem \ref{Th:main23}, we have
$$\bigg|P\left(\sup_{(s,t)\in [0,S]\times[T-\delta_{u},T]}X(s,t)\leq u\right)
-P\left(\sup_{(s,t)\in \cup I_{j}\times [T-\delta_{u},T]}X(s,t)\leq u\right)\bigg|\rightarrow 0,$$
as $u\rightarrow\infty$ and $\epsilon\rightarrow 0$.
\EL

\textbf{Proof:} By applying Theorem \ref{Th:main1} and Lemma \ref{Le:A1}, we have
\begin{eqnarray*}
&&\bigg|P\left(\sup_{(s,t)\in [0,S]\times[T-\delta_{u},T]}X(s,t)\leq u\right)
-P\left(\sup_{(s,t)\in \cup I_{j}\times [T-\delta_{u},T]}X(s,t)\leq u\right)\bigg|\\
&&\leq P\left(\sup_{(s,t)\in (\cup I_{j}^{\epsilon}\cup J)\times [T-\delta_{u},T]}X(s,t)> u\right)\\
&&\leq \sum_{j=1}^{\lfloor S\rfloor}P\left(\sup_{(s,t)\in I_{j}^{\epsilon}\times [T-\delta_{u},T]}X(s,t)> u\right)+
P\left(\sup_{(s,t)\in J\times [T-\delta_{u},T]}X(s,t)> u\right)\\
&&\thicksim \sum_{j=1}^{\lfloor S\rfloor}P\left(\sup_{(s,t)\in I_{j}^{\epsilon}\times [0,T]}X(s,t)> u\right)+
P\left(\sup_{(s,t)\in J\times [0,T]}X(s,t)> u\right)\\
&&\leq (\lfloor S\rfloor\epsilon +1)\mu(u).
\end{eqnarray*}
Noting that by the definitions of $a_{S}$  and $b_{S}$, we have
$S\mu(u)=O(1)$ as  $u\rightarrow\infty$, thus the result follows by letting $\epsilon\rightarrow0$. $\Box$

Let in the following $q_{i}=du^{-2/\alpha_{i}}$ for some $d>0$.

\BL\label{Le:A3} Under the conditions of Theorem \ref{Th:main23}, we have  for any $j=1,2,\cdots,\lfloor S\rfloor$
$$ \bigg|P\left(\sup_{(s,t)\in I_{j}\times [T-\delta_{u},T]}X(s,t)\leq u\right)
-P\left(\sup_{(kq_{1},lq_{2})\in  I_{j}\times [T-\delta_{u},T]}X(kq_{1},lq_{2})\leq u\right)\bigg|\leq K\rho(d)\mu(u)$$
as $u\rightarrow \infty$, where $\rho(d)\rightarrow0$ as $d\rightarrow0$.
\EL

\textbf{Proof:} Without loss of generality, we only show the case $j=1$.\\
 \underline{Case $\beta>\max(\alpha_{1},\alpha_{2})$:}  For simplicity, we only consider the case that $\alpha_{1}=\alpha_{2}=:\alpha$.
Choose first a constant $\alpha_{0}\in (\alpha,\beta)$ and denote that
$$\Delta_{ij}=\Delta_{i}\times\Delta_{j},\ \ \Delta_{ij}^{T}=\Delta_{i}\times(T-\Delta_{j}),\ \ \mbox{with}\ \
\Delta_{i}=[iu^{-\frac{2}{\alpha_{0}}},(i+1)u^{-\frac{2}{\alpha_{0}}}].$$
Set further
$$N_{1}(u)=\lfloor (1-\epsilon)u^{\frac{2}{\alpha_{0}}}\rfloor+1,\ \ N_{2}(u)=\lfloor(\ln u)^{\frac{2}{\beta}}u^{\frac{2}{\alpha_{0}}-\frac{2}{\beta}}\rfloor+1.$$
For any $\varepsilon\in(0,1)$, let $\{\eta_{\pm\varepsilon}(s,t),(s,t)\in[0,\infty)^{2}\}$ be centered homogeneous Gaussian random fields with
covariance functions
$$r_{\pm\varepsilon}(s,t)=\exp\left(-(1\pm \varepsilon)^{\alpha}\big(|a_{1}s|^{\alpha}+|a_{2}t+a_{3}s|^{\alpha}\big)\right),\ \ (s,t)\in[0,\infty)^{2}$$
From the proof of case i) of D\c{e}bicki et al. (2016), it is easy to show that (letting $q=du^{-\frac{2}{\alpha}}$)
\begin{eqnarray}
\label{A1}
&&\sum_{i=0}^{N_{1}(u)}\sum_{j=0}^{N_{2}(u)}P\left(\sup_{(s,t)\in \Delta_{ij}}\eta_{+\varepsilon}(s,T-t)> u_{j-}\right)\nonumber\\
&&\geq \sum_{i=0}^{N_{1}(u)}\sum_{j=0}^{N_{2}(u)}P\left(\sup_{(s,t)\in \Delta_{ij}}\frac{X(s,T-t)}{\sigma(s,T-t)}> u_{j-}\right)\nonumber\\
&&\geq \sum_{i=0}^{N_{1}(u)}\sum_{j=0}^{N_{2}(u)}P\left(\sup_{(s,t)\in \Delta_{ij}^{T}}X(s,t)> u\right)\nonumber\\
&&\geq P\left(\sup_{(s,t)\in I_{1}\times [T-\delta_{u},T]}X(s,t)> u\right)\nonumber\\
&&\geq P\left(\sup_{(kq,lq)\in  I_{1}\times [T-\delta_{u},T]}X(kq,lq)> u\right)\nonumber\\
&&\geq \sum_{i=0}^{N_{1}(u)-1}\sum_{j=0}^{N_{2}(u)-1}P\left(\sup_{(kq,lq)\in \Delta_{ij}^{T}}X(kq,lq)> u\right)-\Sigma_{1}(u)\nonumber\\
&&\geq\sum_{i=0}^{N_{1}(u)-1}\sum_{j=0}^{N_{2}(u)-1}P\left(\sup_{(kq,lq)\in \Delta_{ij}}\frac{X(kq,T-lq)}{\sigma(kq,T-lq)}> u_{j+}\right)-\Sigma_{1}(u)\nonumber\\
&&\geq \sum_{i=0}^{N_{1}(u)-1}\sum_{j=0}^{N_{2}(u)-1}P\left(\sup_{(kq,lq)\in \Delta_{ij}}\eta_{-\varepsilon}(kq,T-lq)> u_{j+}\right)-\Sigma_{1}(u),
\end{eqnarray}
where
$$u_{j-}=u(1+b(1-\varepsilon)(ju^{-\frac{2}{\alpha_{0}}})^{\beta}),\ \ u_{j+}=u(1+b(1+\varepsilon)((j+1)u^{-\frac{2}{\alpha_{0}}})^{\beta}),$$
and
$$\Sigma_{1}(u)=\sum_{0\leq i,i'\leq N_{1}(u)-1,}\sum_{0\leq j,j'\leq N_{2}(u)-1}P\left(\sup_{(s,t)\in \Delta_{ij}^{T}}X(s,t)> u,
\sup_{(s,t)\in \Delta_{i'j'}^{T}}X(s,t)> u\right).$$
We also can get the following results from the above mentioned paper
\begin{eqnarray}
\label{A2}
\Sigma_{1}(u)=o(\mu(u))
\end{eqnarray}
as $u\rightarrow\infty$ and
\begin{eqnarray}
\label{A3}
&&\sum_{i=0}^{N_{1}(u)}\sum_{j=0}^{N_{2}(u)}P\left(\sup_{(s,t)\in \Delta_{ij}}\eta_{+\varepsilon}(s,T-t)> u_{j-}\right)\nonumber\\
&&\ \ \ \sim \sum_{i=0}^{N_{1}(u)-1}\sum_{j=0}^{N_{2}(u)-1}P\left(\sup_{(s,t)\in \Delta_{ij}}\eta_{-\varepsilon}(s,T-t)> u_{j+}\right)\sim \mu(u),
\end{eqnarray}
as $u\rightarrow\infty$ and $\varepsilon\rightarrow 0$.
 For the homogeneous Gaussian random fields $\eta_{\pm\varepsilon}(s,t)$, by Lemma \ref{Le:A62} in the Appendix, we use the following estimate
\begin{eqnarray*}
&&\bigg|P\left(\sup_{(s,t)\in \Delta_{ij}}\eta_{\pm\varepsilon}(s,T-t)> u\right)
-P\left(\sup_{(kq,lq)\in \Delta_{ij}}\eta_{\pm\varepsilon}(kq,T-lq)> u\right)\bigg|\\
&& \leq \rho(d)a_{1}a_{2}u^{\frac{4}{\alpha}-\frac{4}{\alpha_{0}}}\Psi(u)(1+g(u)),
\end{eqnarray*}
where $\rho(d)\rightarrow 0$ as $d\rightarrow 0$ and $g(u)\rightarrow0$ as $u\rightarrow\infty$.
Denote by $G(u)=1+\sup_{v\geq u}|g(u)|\rightarrow1$ as $u\rightarrow\infty$. Then
$u/u_{j\pm}\rightarrow 1$ as $u\rightarrow\infty$ uniformly in $j$ and also
\begin{eqnarray*}
&&\bigg|P\left(\sup_{(s,t)\in \Delta_{ij}}\eta_{\pm\varepsilon}(s,T-t)> u_{j\pm}\right)
-P\left(\sup_{(kq,lq)\in \Delta_{ij}}\eta_{\pm\varepsilon}(kq,T-lq)> u_{j\pm}\right)\bigg|\\
&& \leq \rho(d)a_{1}a_{2}u_{j\pm}^{\frac{4}{\alpha}-\frac{4}{\alpha_{0}}}\Psi(u_{j\pm})G(u).
\end{eqnarray*}
Thus, there exists $K>0$ such that
\begin{eqnarray}
\label{A4}
&&\bigg|\sum_{i=0}^{N_{1}(u)-1}\sum_{j=0}^{N_{2}(u)-1}P\left(\sup_{(s,t)\in \Delta_{ij}}\eta_{-\varepsilon}(s,T-t)> u_{j+}\right)
-\sum_{i=0}^{N_{1}(u)-1}\sum_{j=0}^{N_{2}(u)-1}P\left(\sup_{(kq,lq)\in \Delta_{ij}}\eta_{-\varepsilon}(kq,T-lq)> u_{j+}\right)\bigg|\nonumber\\
&& \leq K\rho(d)u^{\frac{4}{\alpha}-\frac{2}{\beta}}\Psi(u).
\end{eqnarray}
Now it follows from (\ref{A1}-\ref{A4}) that
\begin{eqnarray*}
&&\bigg| P\left(\sup_{(s,t)\in I_{1}\times[T-\delta_{u},T]}X(s,t)> u\right)-
 P\left(\sup_{(kq,lq)\in I_{1}\times[T-\delta_{u},T]}X(kq,lq)> u\right)\bigg|\\
 && \leq\sum_{i=0}^{N_{1}(u)}\sum_{j=0}^{N_{2}(u)}P\left(\sup_{(s,t)\in \Delta_{ij}}\eta_{+\varepsilon}(s,T-t)> u_{j-}\right)-\sum_{i=0}^{N_{1}(u)-1}\sum_{j=0}^{N_{2}(u)-1}P\left(\sup_{(kq,lq)\in \Delta_{ij}}\eta_{-\varepsilon}(kq,T-lq)> u_{j+}\right)+\Sigma_{1}(u)\\
 && \sim \sum_{i=0}^{N_{1}(u)-1}\sum_{j=0}^{N_{2}(u)-1}P\left(\sup_{(s,t)\in \Delta_{ij}}\eta_{-\varepsilon}(s,T-t)> u_{j+}\right)-
 \sum_{i=0}^{N_{1}(u)-1}\sum_{j=0}^{N_{2}(u)-1}P\left(\sup_{(kq,lq)\in \Delta_{ij}}\eta_{-\varepsilon}(kq,T-lq)> u_{j+}\right)\\
 &&\leq K\rho(d)u^{\frac{4}{\alpha}-\frac{2}{\beta}}\Psi(u).
\end{eqnarray*}

\underline{Case $\beta=\alpha_{1}=\alpha_{2}$:} For simplicity, set $\alpha=\alpha_{1}=\alpha_{2}$. Let $S_{0},T_{0}$ be two positive constants
and define
$$\widehat{\Delta}_{i}=[iS_{0}u^{-\frac{2}{\alpha}},(i+1)S_{0}u^{-\frac{2}{\alpha}}],\ \ i=0,1,\cdots,\widehat{N}_{1}(u),\ \
\widetilde{\Delta}_{j}=[jT_{0}u^{-\frac{2}{\alpha}},(j+1)T_{0}u^{-\frac{2}{\alpha}}],\ \ j=0,1,\cdots,\widetilde{N}_{2}(u),$$
$$\overline{\Delta}_{ij}=\widehat{\Delta}_{i}\times\widetilde{\Delta}_{j},
\ \ \overline{\Delta}_{ij}^{T}=\widehat{\Delta}_{i}\times(T-\widetilde{\Delta}_{j}),$$
where
$$\widehat{N}_{1}(u)=\lfloor\frac{1-\epsilon}{S_{0}}u^{\frac{2}{\alpha}}\rfloor+1,\ \ \widetilde{N}_{2}(u)=\lfloor\frac{(\ln u)^{\frac{2}{\beta}}}{T_{0}}u^{\frac{2}{\alpha}}\rfloor+1.$$
From the proof of case ii) of D\c{e}bicki et al. (2016) again, it is easy to show that (letting $q=du^{-\frac{2}{\alpha}}$)
\begin{eqnarray}
\label{B1}
&&\Sigma_{2}(u)+\sum_{i=0}^{\widehat{N}_{1}(u)}P\left(\sup_{(s,t)\in \overline{\Delta}_{i0}^{T}}X(s,t)> u\right)\nonumber\\
&&\geq P\left(\sup_{(s,t)\in I_{1}\times[T-\delta_{u},T]}X(s,t)> u\right)\nonumber\\
&&\geq P\left(\sup_{(kq,lq)\in I_{1}\times[T-\delta_{u},T]}X(kq,lq)> u\right)\nonumber\\
&&\geq \sum_{i=0}^{\widehat{N}_{1}(u)-1}P\left(\sup_{(kq,lq)\in \overline{\Delta}_{i0}^{T}}X(kq,lq)> u\right)-\Sigma_{3}(u),
\end{eqnarray}
where
$$\Sigma_{2}(u)=\sum_{i=0}^{\widehat{N}_{1}(u)}\sum_{j=1}^{\widetilde{N}_{2}(u)}P\left(\sup_{(s,t)\in \overline{\Delta}_{ij}^{T}}X(s,t)> u\right)=o(\mu(u)),$$
$$\Sigma_{3}(u)=\sum_{0<i<i'<\widehat{N}_{1}(u)-1}P\left(\sup_{(s,t)\in \overline{\Delta}_{i0}^{T}}X(s,t)> u,
\sup_{(s,t)\in \overline{\Delta}_{i'0}^{T}}X(s,t)> u\right)=o(\mu(u)),$$
as $u\rightarrow\infty$. We can also get the following results by Lemma 2.1 of D\c{e}bicki et al. (2016)
$$P\left(\sup_{(s,t)\in \overline{\Delta}_{i0}^{T}}X(s,t)> u\right)\thicksim P\left(\sup_{(s,t)\in \overline{\Delta}_{i0}}\frac{\widetilde{\eta}(s,t)}{1+bt^{\beta}}> u\right)\thicksim \mathcal{H}_{Y_{1}}^{b}[S_{0},T_{0}]\Psi(u)$$
and
$$\sum_{i=0}^{\widehat{N}_{1}(u)}P\left(\sup_{(s,t)\in \overline{\Delta}_{i0}^{T}}X(s,t)> u\right)\thicksim \sum_{i=0}^{\widehat{N}_{1}(u)-1}P\left(\sup_{(s,t)\in \overline{\Delta}_{i0}^{T}}X(s,t)> u\right)\thicksim
\frac{(1-\epsilon)}{S_{0}}u^{\frac{2}{\alpha}}\mathcal{H}_{Y_{1}}^{b}[S_{0},T_{0}]\Psi(u)$$
as $u\rightarrow\infty$,
where $\{\widetilde{\eta}(s,t),(s,t)\in[0,\infty)^{2}\}$ is a  centered homogeneous Gaussian random fields with
covariance functions
$$r(s,t)=\exp\left(-\big(|a_{1}s|^{\alpha}+|a_{2}t-a_{3}s|^{\alpha}\big)\right),\ \ (s,t)\in[0,\infty)^{2}.$$
Since $X(s,t)$ is homogeneous with respect
to $s$, we have
\begin{eqnarray}
\label{B2}
0&\leq& P\left(\sup_{(s,t)\in  I_{1}\times[T-\delta_{u},T]}X(s,t)> u\right)
- P\left(\sup_{(kq,lq)\in I_{1}\times[T-\delta_{u},T]}X(kq,lq)> u\right)\nonumber\\
&\leq &\Sigma_{2}(u)+\sum_{i=0}^{\widehat{N}_{1}(u)}P\left(\sup_{(s,t)\in \overline{\Delta}_{i0}^{T}}X(s,t)> u\right)
- \sum_{i=0}^{\widehat{N}_{1}(u)-1}P\left(\sup_{(kq,lq)\in \overline{\Delta}_{i0}^{T}}X(kq,lq)> u\right)+\Sigma_{3}(u)\nonumber\\
&=&\Sigma_{2}(u)+\widehat{N}_{1}(u)P\left(\sup_{(s,t)\in \overline{\Delta}_{00}^{T}}X(s,t)> u\right)
- (\widehat{N}_{1}(u)-1)P\left(\sup_{(kq,lq)\in \overline{\Delta}_{00}^{T}}X(kq,lq)> u\right)+\Sigma_{3}(u)\nonumber\\
&\leq&(\widehat{N}_{1}(u)-1)P\left(\sup_{(s,t)\in \overline{\Delta}_{00}}\frac{\widetilde{\eta}(s,t)}{1+bt^{\beta}}> u, \sup_{(kq,lq)\in \overline{\Delta}_{00}}\frac{\widetilde{\eta}(kq,lq)}{1+b(lq)^{\beta}}\leq u\right)+\mathcal{H}_{Y_{1}}^{b}[S_{0},T_{0}]\Psi(u)+o(\mu(u)).
\end{eqnarray}
For the constants $a_{1}>0, a_{2}>0, a_{3}\neq 0, b>0$, let (as in Section 2)
$$Y(s,t)=\widetilde{B}_{\alpha}(a_{1}s)+B_{\alpha}(a_{2}t-a_{3}s),\ \ \sigma_{Y}^{2}(s,t)=Var(Y(s,t))$$
and
$$\mathcal{H}_{Y}^{b}[\lambda_{1},\lambda_{2}](d)=E\exp\left(\max_{(kd,ld)\in[0,d\lambda_{1}]\times[0,d\lambda_{2}]}\sqrt{2}Y(kd,ld)-\sigma_{Y}^{2}(kd,ld)-b|ld|^{\beta}\right)\in (0,\infty),$$
where $\widetilde{B}_{\alpha}$ and $B_{\alpha}$ are two independent fBms.
By the same arguments as in the proof of Lemma 6.1 of  D\c{e}bicki et al. (2016), we can show
$$\mathcal{M}_{Y,\alpha}^{b}(d)=\lim_{\lambda_{1}\rightarrow\infty}\lim_{\lambda_{2}\rightarrow\infty}\frac{1}{d\lambda_{1}}\mathcal{H}_{Y}^{b}[\lambda_{1},\lambda_{2}](d)\in (0,\infty).$$
Following the arguments of Lemma 12.2.7 of Leadbetter et al. (1983), we can show that $\lim_{d\rightarrow0}\mathcal{M}_{Y,\alpha}^{b}(d)=\mathcal{M}_{Y,\alpha}^{b}$.
Now, following the arguments of Lemma 6.1 of  D\c{e}bicki et al. (2016)  (see also the proof of Lemma 6.1 of Piterbarg (1996)), we have
\begin{eqnarray}
\label{B3}
&&P\left(\sup_{(s,t)\in \overline{\Delta}_{00}}\frac{\widetilde{\eta}(s,t)}{1+bt^{\beta}}> u, \sup_{(kq,lq)\in \overline{\Delta}_{00}}\frac{\widetilde{\eta}(kq,lq)}{1+b(lq)^{\beta}}\leq u\right)\nonumber\\
&&=\Psi(u)\int_{0}^{+\infty}e^{w}P\bigg(\sup_{(s,t)\in [0,S_{0}]\times[0,T_{0}]}[\sqrt{2}Y(s,t)-\sigma_{Y}^{2}(s,t)-b|t|^{\beta}]> w,\nonumber\\
&&\ \ \ \ \ \ \ \ \ \ \ \
\sup_{(kq,lq)\in [0,S_{0}]\times[0,T_{0}]}[\sqrt{2}Y(kq,lq)-\sigma_{Y}^{2}(kd,ld)-b|lq|^{\beta}]\leq w\bigg)dw (1+o(1))\nonumber\\
&&=\Psi(u)\left(\mathcal{H}_{Y}^{b}[S_{0}/d,T_{0}/d](d)-\mathcal{H}_{Y}^{b}[S_{0},T_{0}]\right)(1+o(1)),
\end{eqnarray}
as $u\rightarrow\infty$.
\COM{where $\chi(s,t)$ is a Gaussian random field with continuous trajectories which is defined by
$$\chi(s,t)=B_{\alpha}(a_{1}s)+\widetilde{B}_{\alpha}(a_{2}t-a_{3}s)$$
with $B_{\alpha}$ and $\widetilde{B}_{\alpha}$ two independent fBm's defined on $\mathbb{R}$ with Hurst
index $\alpha/2\in(0,1]$. Since the sample paths of $\chi(s,t)$ are continuous,
then for a fixed $S_{0},T_{0}$ the probability in (\ref{B3}) tends to zero as $d\rightarrow 0$.
Let's denote by $\rho(d)$ for the integral of (\ref{B3}), then by the dominated convergence theorem we find that
$\rho(d)\rightarrow0$ as $d\rightarrow 0$.} Now, we can conclude that
\begin{eqnarray}
\label{B4}
0&\leq& P\left(\sup_{(s,t)\in I_{1}\times[T-\delta_{u},T]}X(s,t)> u\right)
- P\left(\sup_{(kq,lq)\in I_{1}\times[T-\delta_{u},T]}X(kq,lq)> u\right)\nonumber\\
&\leq& (1-\epsilon)u^{2/\alpha}\Psi(u)\left(\frac{\mathcal{H}_{Y}^{b}[S_{0}/d,T_{0}/d](d)}{S_{0}}-\frac{\mathcal{H}_{Y}^{b}[S_{0},T_{0}]}{S_{0}}\right)+\mathcal{H}_{Y_{1}}^{b}[S_{0},T_{0}]\Psi(u)+o(\mu(u))\nonumber\\
&\leq& K\left(\mathcal{M}_{Y,\alpha}^{b}(d)-\mathcal{M}_{Y,\alpha}^{b}\right)\mu(u)\nonumber\\
&=:&K\rho(d)\mu(u),
\end{eqnarray}
where $\rho(d)\rightarrow0$ as $d\rightarrow0$.

\underline{Case $\beta=\alpha_{2}>\alpha_{1}$:} This case can be proved as case ii) by some obvious changes as follows.
Let $S_{0},T_{0}$ be two positive constants
and define
$$\widehat{\Delta}_{i}=[iS_{0}u^{-\frac{2}{\alpha_{1}}},(i+1)S_{0}u^{-\frac{2}{\alpha_{1}}}],\ \ i=0,1,\cdots,\widehat{N}_{1}(u),\ \
\widetilde{\Delta}_{j}=[jT_{0}u^{-\frac{2}{\alpha_{2}}},(j+1)T_{0}u^{-\frac{2}{\alpha_{2}}}],\ \ j=0,1,\cdots,\widetilde{N}_{2}(u),$$
$$\overline{\Delta}_{ij}=\widehat{\Delta}_{i}\times\widetilde{\Delta}_{j},
\ \ \overline{\Delta}_{ij}^{T}=\widehat{\Delta}_{i}\times(T-\widetilde{\Delta}_{j}),$$
where
$$\widehat{N}_{1}(u)=\lfloor\frac{1-\epsilon}{S_{0}}u^{\frac{2}{\alpha_{1}}}\rfloor+1,\ \ \widetilde{N}_{2}(u)=\lfloor\frac{(\ln u)^{\frac{2}{\beta}}}{T_{0}}u^{\frac{2}{\alpha_{2}}}\rfloor+1.$$
Let $q_{1}=du^{-\frac{2}{\alpha_{1}}}, q_{2}=du^{-\frac{2}{\alpha_{2}}}$, then repeating the proof of case ii) by replacing $kq$ and $lq$
by $kq_{1}$ and $lq_{2}$, we get the desired result.

\underline{Case $\beta<\alpha_{2}=\alpha_{1}$:} For simplicity let $\alpha:=\alpha_{2}=\alpha_{1}$ and $q=du^{-\frac{2}{\alpha}}$. Let's consider the Gaussian process $X(s,T), s\geq 0$. It is easy to check that
$X(s,T), s\geq 0$ is standard stationary Gaussian process, i.e., with mean 0, variance 1. For the covariance function of $X(s,T), s\geq 0$, it holds that
$$r(s,T,s',T)=1-(a_{1}^{\alpha}+|a_{3}|^{\alpha})|s-s'|^{\alpha}(1+o(1))$$
uniformly with respect to $s,s'\in [0,S_{0}]$, as $|s-s'|\rightarrow 0$.
For some constant $a>0$, let
$$\mathcal{H}_{\alpha}^{a}[0,\lambda]=E\exp\left(\max_{ak\in[0,a\lambda]}\sqrt{2}B_{\alpha}(ak)-(ak)^{\alpha}\right)$$
and define
$$\mathcal{H}_{\alpha}(a)=\lim_{\lambda\rightarrow\infty}\frac{\mathcal{H}_{\alpha}^{a}[0,\lambda]}{a\lambda}\in (0,+\infty).$$
Note that $\lim_{a\rightarrow0}\mathcal{H}_{\alpha}(a)=\mathcal{H}_{\alpha}$, see e.g. Leadbetter et al. (1983).
So by Lemmas \ref{Le:A61} and \ref{Le:A62} in the Appendix (for the one dimensional case), we have
\begin{eqnarray}
\label{D1}
P\left(\sup_{s\in  I_{1}}X(s,T)> u\right)
=(1-\epsilon)(a_{1}^{\alpha}+|a_{3}|^{\alpha})^{\frac{1}{\alpha}}\mathcal{H}_{\alpha}u^{\frac{2}{\alpha}}\Psi(u)(1+o(1)),
\end{eqnarray}
\begin{eqnarray}
\label{D2}
P\left(\sup_{kq\in  I_{1}}X(kq,T)> u\right)
=(1-\epsilon)(a_{1}^{\alpha}+|a_{3}|^{\alpha})^{\frac{1}{\alpha}}\mathcal{H}_{\alpha}(d)u^{\frac{2}{\alpha}}\Psi(u)(1+o(1))
\end{eqnarray}
and
\begin{eqnarray}
\label{D3}
 \bigg|P\left(\sup_{s\in  I_{1}}X(s,T)> u\right)
-P\left(\sup_{kq\in  I_{1}}X(kq,T)> u\right)\bigg|\leq K\rho(d)u^{\frac{2}{\alpha}}\Psi(u),
\end{eqnarray}
as $u\rightarrow \infty$, where $\rho(d)=\mathcal{H}_{\alpha}(d)-\mathcal{H}_{\alpha}$.
By repeating the proof of iv) of D\c{e}bicki et al. (2016), it is easy to show that
\begin{eqnarray}
\label{D4}
P\left(\sup_{(kq,lq)\in  I_{1}\times[T-\delta_{u},T]}X(kq,lq)> u\right)=(1-\epsilon)(a_{1}^{\alpha}+|a_{3}|^{\alpha})^{\frac{1}{\alpha}}\mathcal{H}_{\alpha}(d)u^{\frac{2}{\alpha}}\Psi(u)(1+o(1)).
\end{eqnarray}
Write
\begin{eqnarray*}
\label{D5}
&&\bigg|P\left(\sup_{(s,t)\in I_{1}\times[T-\delta_{u},T]}X(s,t)> u\right)
-P\left(\sup_{(kq,lq)\in I_{1}\times[T-\delta_{u},T]}X(kq,lq)> u\right)\bigg|\nonumber\\
&&\leq \bigg|P\left(\sup_{(s,t)\in I_{1}\times[T-\delta_{u},T]}X(s,t)> u\right)
-P\left(\sup_{s\in I_{1}}X(s,T)> u\right)\bigg|\nonumber\\
&&+\bigg|P\left(\sup_{s\in I_{1}}X(s,T)> u\right)-P\left(\sup_{kq\in I_{1}}X(kq,T)> u\right)\bigg|\nonumber\\
&&+\bigg|P\left(\sup_{kq\in I_{1}}X(kq,T)> u\right)-P\left(\sup_{(kq,lq)\in I_{1}\times[T-\delta_{u},T]}X(kq_{1},lq_{2})> u\right)\bigg|\nonumber\\
&&=:M_{1}+M_{2}+M_{3},
\end{eqnarray*}
where $M_{1}=o(\mu(u))$ by iv) of Theorem \ref{Th:main1} and (\ref{D1}),
$M_{2}=K\rho(d)u^{\frac{2}{\alpha}}\Psi(u)$  by (\ref{D3}) and $M_{3}=o(\mu(u))$ by (\ref{D2}) and (\ref{D4}) as $u\rightarrow\infty$.

\underline{Case $\beta<\alpha_{2}$ and $\alpha_{1}<\alpha_{2}$:} The proof is the same as that of \underline{Case $\beta<\alpha_{2}=\alpha_{1}$}.

\underline{Case $\beta=\alpha_{1}>\alpha_{2}$ and case $\beta<\alpha_{1}$ and $\alpha_{2}<\alpha_{1}$:} These two cases can be proved by the same arguments as for the third and fifth cases after some time scaling as in
D\c{e}bicki et al. (2016), so we omit the details. $\Box$

\BL\label{Le:A4} Under the conditions of Theorem \ref{Th:main23}, we have
\begin{eqnarray}
\bigg|P\left(\sup_{(s,t)\in  \cup I_{j}\times[T-\delta_{u},T]}X(s,t)\leq u\right)
-P\left(\sup_{(kq_{1},lq_{2})\in \cup I_{j}\times[T-\delta_{u},T]}X(kq_{1},lq_{2})\leq u\right)\bigg|\leq K\rho(d)S\mu(u)
\end{eqnarray}
as $u\rightarrow \infty$.
\EL

\textbf{Proof:}  By Lemma \ref{Le:A3}, we have
\begin{eqnarray*}
&&\bigg|P\left(\sup_{(s,t)\in  \cup I_{j}\times[T-\delta_{u},T]}X(s,t)\leq u\right)
-P\left(\sup_{(kq_{1},lq_{2})\in \cup I_{j}\times[T-\delta_{u},T]}X(kq_{1},lq_{2})\leq u\right)\bigg|\\
&& \leq S \max_{j}\bigg|P\left(\sup_{(s,t)\in  I_{j}\times[T-\delta_{u},T]}X(s,t)\leq u\right)
-P\left(\sup_{(kq_{1},lq_{2})\in  I_{j}\times[T-\delta_{u},T]}X(kq_{1},lq_{2})\leq u\right)\bigg|\\
&&\leq K\rho(d)S\mu(u),
\end{eqnarray*}
which completes the proof. $\Box$

\BL\label{Le:A5} Under the conditions of Theorem \ref{Th:main23}, we have
$$ \bigg|P\left(\sup_{(kq_{1},lq_{2})\in \cup I_{j}\times[T-\delta_{u},T]}X(kq_{1},lq_{2})\leq u\right)
-\prod_{j=1}^{\lfloor S\rfloor}P\left(\sup_{(kq_{1},lq_{2})\in  I_{j}\times[T-\delta_{u},T]}X(kq_{1},lq_{2})\leq u\right)\bigg|\rightarrow0,$$
as $u\rightarrow \infty$.
\EL

\textbf{Proof:}
Applying Berman's inequality (see e.g.
Piterbarg (1996)) we have
\begin{eqnarray*}
&&\bigg|P\left(\sup_{(kq_{1},lq_{2})\in \cup I_{j}\times[T-\delta_{u},T]}X(kq_{1},lq_{2})\leq u\right)
-\prod_{j=1}^{\lfloor S\rfloor}P\left(\sup_{(kq_{1},lq_{2})\in  I_{j}\times[T-\delta_{u},T]}X(kq_{1},lq_{2})\leq u\right)\bigg|\\
&&=\bigg|P\left(\sup_{(kq_{1},lq_{2})\in \cup I_{j}\times[T-\delta_{u},T]}\frac{X(kq_{1},lq_{2})}{\sigma(lq_{2})}\leq \frac{u}{\sigma(lq_{2})}\right)
-\prod_{j=1}^{\lfloor S\rfloor}P\left(\sup_{(kq_{1},lq_{2})\in I_{j}\times[T-\delta_{u},T]}\frac{X(kq_{1},lq_{2})}{\sigma(lq_{2})}\leq \frac{u}{\sigma(lq_{2})}\right)\bigg|\\
&&\leq \sum_{j\neq j'}\sum_{(kq_{1},lq_{2})\in I_{j}\times[T-\delta_{u},T]\atop (k'q_{1},l'q_{2})\in  I_{j'}\times[T-\delta_{u},T]}
|r(kq_{1},lq_{2},k'q_{1},l'q_{2})|\exp\left(-\frac{(\sigma^{-2}(lq_{2})+\sigma^{-2}(l'q_{2}))u^{2}}{2(1+r(kq_{1},lq_{2},k'q_{1},l'q_{2}))}\right)\\
&&\leq \sum_{j\neq j'}\sum_{(kq_{1},lq_{2})\in I_{j}\times[T-\delta_{u},T]\atop (k'q_{1},l'q_{2})\in  I_{j'}\times[T-\delta_{u},T]}
|r(kq_{1},lq_{2},k'q_{1},l'q_{2})|\exp\left(-\frac{u^{2}}{1+r(kq_{1},lq_{2},k'q_{1},l'q_{2})}\right).
\end{eqnarray*}
Since $|kq_{1}-k'q_{1}|\geq \epsilon$ by definition,
$r(kq_{1},lq_{2},k'q_{1},l'q_{2})\leq\delta<1$. Set
$\gamma<(1-\delta)/(1+\delta)$ and split the last sum into two parts
$W_{1}$ and $W_{2}$ with $|kq_{1}-k'q_{1}|< S^{\gamma}$ and
$|kq_{1}-k'q_{1}|\geq S^{\gamma}$, respectively. For the first
sum there are $S^{1+\gamma}/q_{2}^{2}$ combinations of two points
$kq_{1}, k'q_{1}\in \cup_{j} \mathbf{I}_{j}$. Together with the $lq_{2}$
combinations there are $(S^{1+\gamma}/q_{1}^{2})(\delta_{u}/q_{2}^{2})$ terms in
the sum $W_{1}$. Note that
$$S\mu(u)=O(1),\ \ u\rightarrow\infty,$$
which implies for case i)
$$u^{2}=2\ln S+(\frac{2}{\alpha_{1}}+\frac{2}{\alpha_{2}}-\frac{2}{\beta}-1)\ln\ln S+O(1);$$
for case ii)-v)
$$u^{2}=2\ln S+(\frac{2}{\alpha_{1}}-1)\ln\ln S+O(1);$$
for case vi)-vii)
$$u^{2}=2\ln S+(\frac{2}{\alpha_{2}}-1)\ln\ln S+O(1).$$

Thus, $W_{1}$ is bounded by
\begin{eqnarray*}
&&\delta \frac{S^{1+\gamma}\delta^{2}(u)}{q_{1}^{2}q_{2}^{2}}\exp\left(-\frac{u^{2}}{1+\delta}\right) \\
&&\leq\delta\exp\left((1+\gamma)\ln S+(\frac{1}{\alpha_{1}}+\frac{1}{\alpha_{2}})\ln\ln S-\frac{2(1+o(1))}{1+\delta}\ln S\right)\\
&&=\delta\exp\left((\ln S)\left[(1+\gamma)-\frac{2(1+o(1))}{1+\delta}+\frac{(\frac{1}{\alpha_{1}}+\frac{1}{\alpha_{2}})\ln\ln S}{\ln S}\right]\right)\rightarrow0
\end{eqnarray*}
 \wE{as $S \to \IF$} since
$1+\gamma<2/(1+\delta)$ by the choice of
$\gamma$.\\
For the second sum $W_{2}$ with $|kq_{1}-k'q_{1}|\geq S^{\gamma}$, we use that
$$\sup_{|kq_{1}-k'q_{1}|\geq S^{\gamma}}r(kq_{1},lq_{2},k'q_{1},l'q_{2})(\ln S)^{c}=o(1),$$
as $S\rightarrow\infty$.
 In this case there $(S/q_{1})^{2}$ many
combinations of two points $kq_{1}, k'q_{1}\in \cup_{i} \mathbf{I}_{i}$.
Hence $W_{2}$ is bounded by
\begin{eqnarray*}
&&R(S):=\frac{o(1)}{(\ln S)^{c}}\frac{S^{2}}{q_{1}^{2}}\frac{\delta^{2}(u)}{q_{2}^{2}}\exp\left(-\frac{u^{2}}{1+o(1)/\ln S}\right)\\
&&\leq  C\exp\left(2\ln S+(\frac{2}{\alpha_{1}}+\frac{2}{\alpha_{2}}-\frac{2}{\beta}-1)\ln\ln S+\frac{4}{\beta}\ln\ln\ln S-(c-1)\ln\ln S-\frac{u^{2}}{1+o(1)/\ln S}\right).
\end{eqnarray*}
For case i), by assumption {\bf  A4}, $c>1$, we have
\begin{eqnarray*}
&&R(S)\leq  C\exp\bigg(2\ln S+(\frac{2}{\alpha_{1}}+\frac{2}{\alpha_{2}}-\frac{2}{\beta}-1)\ln\ln S+\frac{4}{\beta}\ln\ln\ln S-(c-1)\ln\ln S\\
&&\ \ \ \ \ \ \ \ \ \ \ \ \ \ \ \ \ \ \ \ \ -\frac{(1+o(1))}{1+o(1)/\ln S}[2\ln S+(\frac{2}{\alpha_{1}}+\frac{2}{\alpha_{2}}-\frac{2}{\beta}-1)\ln\ln S]\bigg)\\
&& \leq  C\exp\left(-(c-1)\ln\ln S+\frac{4}{\beta}\ln\ln\ln S+o(1)\right)
\rightarrow0,
\end{eqnarray*}
as $S\rightarrow\infty$, since $c>1$.
For cases ii)-iii), noting that $c>1$, we have
\begin{eqnarray*}
&&R(S)\leq  C\exp\bigg(2\ln S+(\frac{2}{\alpha_{1}}+\frac{2}{\alpha_{2}}-\frac{2}{\beta}-1)\ln\ln S+\frac{4}{\beta}\ln\ln\ln S-(c-1)\ln\ln S\\
&&\ \ \ \ \ \ \ \ \ \ \ \ \ \ \ \ \ \ \ \ \ -\frac{(1+o(1))}{1+o(1)/\ln S}[2\ln S+(\frac{2}{\alpha_{1}}-1)\ln\ln S]\bigg)\\
&& \leq  C\exp\left((\frac{2}{\alpha_{2}}-\frac{2}{\beta})\ln\ln S-(c-1)\ln\ln S+\frac{4}{\beta}\ln\ln\ln S+o(1)\right)
\rightarrow0,
\end{eqnarray*}
as $S\rightarrow\infty$, since $\beta= \alpha_{2}$ and $c>1$.
For cases iv)-v), noting that $c=1$, we have
\begin{eqnarray*}
&&R(S)\leq  C\exp\bigg(2\ln S+(\frac{2}{\alpha_{1}}+\frac{2}{\alpha_{2}}-\frac{2}{\beta}-1)\ln\ln S+\frac{4}{\beta}\ln\ln\ln S\\
&&\ \ \ \ \ \ \ \ \ \ \ \ \ \ \ \ \ \ \ \ \ -\frac{(1+o(1))}{1+o(1)/\ln S}[2\ln S+(\frac{2}{\alpha_{1}}-1)\ln\ln S]\bigg)\\
&& \leq  C\exp\left((\frac{2}{\alpha_{2}}-\frac{2}{\beta})\ln\ln S+\frac{4}{\beta}\ln\ln\ln S+o(1)\right)
\rightarrow0,
\end{eqnarray*}
as $S\rightarrow\infty$, since $\beta< \alpha_{2}$.
For cases vi), we have
\begin{eqnarray*}
&&R(S)\leq  C\exp\bigg(2\ln S+(\frac{2}{\alpha_{1}}+\frac{2}{\alpha_{2}}-\frac{2}{\beta}-1)\ln\ln S+\frac{4}{\beta}\ln\ln\ln S-(c-1)\ln\ln S\\
&&\ \ \ \ \ \ \ \ \ \ \ \ \ \ \ \ \ \ \ \ \ -\frac{(1+o(1))}{1+o(1)/\ln S}[2\ln S+(\frac{2}{\alpha_{2}}-1)\ln\ln S]\bigg)\\
&& \leq  C\exp\left((\frac{2}{\alpha_{1}}-\frac{2}{\beta})\ln\ln S-(c-1)\ln\ln S+\frac{4}{\beta}\ln\ln\ln S+o(1)\right)
\rightarrow0,
\end{eqnarray*}
as $S\rightarrow\infty$, since $\beta= \alpha_{1}$ and $c>1$.
For cases vii), noting that $c=1$, we have
\begin{eqnarray*}
&&R(S)\leq  C\exp\bigg(2\ln S+(\frac{2}{\alpha_{1}}+\frac{2}{\alpha_{2}}-\frac{2}{\beta}-1)\ln\ln S+\frac{4}{\beta}\ln\ln\ln S\\
&&\ \ \ \ \ \ \ \ \ \ \ \ \ \ \ \ \ \ \ \ \ -\frac{(1+o(1))}{1+o(1)/\ln S}[2\ln S+(\frac{2}{\alpha_{2}}-1)\ln\ln S]\bigg)\\
&& \leq  C\exp\left((\frac{2}{\alpha_{1}}-\frac{2}{\beta})\ln\ln S+\frac{4}{\beta}\ln\ln\ln S+o(1)\right)
\rightarrow0,
\end{eqnarray*}
as $S\rightarrow\infty$, since $\beta< \alpha_{1}$. $\Box$

\textbf{Proof of Theorem \ref{Th:main23}:} Recall that $u=u(x)=a_{S}^{-1}x+b_{S}$.  By the stationarity of $X(s,t)$ with respect to the first
component, Lemma \ref{Le:A1}, Theorem \ref{Th:main1} and the choice of $a_{S}, b_{S}$, we have
\BQNY
\prod_{j=1}^{\lfloor S\rfloor}P\left(\max_{(s,t)\in I_{j}\times[T-\delta_{u},T]}X(s,t) \leq
u\right)&\sim&
\exp\left(-\lfloor S\rfloor P\left(\max_{(s,t)\in I_{1}\times[T-\delta_{u},T]} X(s,t) > u\right)\right)\\
&\sim&\exp\left(-\lfloor S\rfloor(1-\epsilon)\mu(u)\right)\\
&\rightarrow& \exp(-e^{-x}), \quad \epsilon\downarrow0,\ \ S \to \IF.
\EQNY
Further, by Lemmas \ref{Le:A1}-\ref{Le:A5}, it holds that as $S \to \IF$
\begin{eqnarray*}
P\left( \max_{(s,t)\in[0,S]\times[0,T]} X(s,t) \leq u\right)
&\sim&P\left(\max_{(s,t)\in\cup_{j} I_{j}\times[T-\delta_{u},T]} X(s,t) \leq u\right)\\
&\sim&P\left(\max_{(kq_{1},lq_{2})\in\cup_{j} I_{j}}X(kq_{1},lq_{2})\leq u\right)\\
&\sim&\prod_{j=1}^{\lfloor S\rfloor}P\left(\max_{(kq_{1},lq_{2})\in I_{j}\times[T-\delta_{u},T]} X(kq_{1},lq_{2})\leq u\right).
\end{eqnarray*}
Therefore, the claim
follows. $\Box$

\textbf{Proof of Proposition \ref{Pro:main3:1}:}
In the paper of D\c{e}bicki et al. (2016), it is shown that
the standard deviation function of $Z$ satisfies assumption {\bf A1} and
the correlation function  of $Z$  satisfies assumption {\bf A2}. It is also shown that
assumption {\bf A3} holds for $Z$. So, in order to prove this proposition, it suffices to show
assumption {\bf A4} holds.
For the correlation function $r_{Z}(s,t,s',t')$ of $Z$, we have
$$r_{Z}(s,t,s',t')=r_{X}(|s+t-s'-t'|)-r_{X}(|s-s'-t'|)-r_{X}(|s+t-s'|)+r_{X}(|s-s'|).$$
Since $r_{X}(t)$ is twice continuously differentiable in $(0,\infty)$, we have
$$|r_{X}(|s+t-s'-t'|)-r_{X}(|s-s'-t'|)-r_{X}(|s+t-s'|)+r_{X}(|s-s'|)|\leq C \ddot{r}_{X}(s-s')$$
for $t,t'\in[0,T]$ as $s-s'\rightarrow\infty$.
Now using the condition that $\ddot{r}_{X}(t)(\ln t)^{c}\rightarrow 0$ as $t\rightarrow \infty$, we show that
assumption {\bf A4} holds. $\Box$

\textbf{Proof of Proposition \ref{pro:main3:1}:} We check that assumptions ${\bf A1-A4}$ hold. Using the stationarity of the increments of $X(t)$ and {\bf C1}, it follows that the variance $\sigma_{Z}^{2}(s,t)$ of $Z(s,t)$ attains its maximum on $[0,T]$  at the unique point $T$, and further
$$\sigma_{Z}(s,t)=\sigma_{X}(t)=1-b(T-t)^{\beta}(1+o(1)),\ \ t\uparrow T$$
holds for some $\beta, b>0$.\\
Notice that for the process $X(t)$ with stationary increments
$$Cov(X(t),X(s))=\frac{1}{2}[\sigma^{2}_{X}(t)+\sigma^{2}_{X}(s)-\sigma^{2}_{X}(|t-s|)].$$
Thus, using the stationarity of the increments of $X(t)$ again, we have for correlation function of $Z(s,t)$
$$r_{Z}(s,t,s',t')=\frac{1}{2\sigma_{X}(t)\sigma_{X}(t')}[-\sigma_{X}^{2}(|s+t-s'-t'|)+\sigma_{X}^{2}(|s-s'-t'|)+\sigma_{X}^{2}(|s-s'+t'|)-\sigma_{X}^{2}(|s-s'|)].$$
It follows from  {\bf C2} that
$$r_{Z}(s,t,s',t')=1-\frac{1}{2}[(a|s+t-s'-t'|)^{\alpha}+(a|s-s'|)^{\alpha}](1+o(1)),$$
as $t,t\rightarrow T$ and $|s-s'|\rightarrow 0$.  {\bf A3} holds obviously. Thus, by Theorem \ref{Th:main1}, the first assertion of Proposition \ref{pro:main3:1} holds.\\
By Taylor expansions, it is straightforward to verify that
$$|r_{Z}(s,t,s',t')|\leq C\ddot{\sigma}_{X}^{2}(|s-s'|)$$
as $|s-s'|\rightarrow\infty$, which combined with {\bf C3} implies {\bf A4}. Thus, by Theorem \ref{Th:main23}, the second assertion holds.
$\Box$

\section{Appendix}

 Let $\{\xi(\mathbf{t}):\mathbf{t}\geq
\mathbf{0}\}$ denote a  two dimensional homogeneous Gaussian field with covariance
function
$$r_{\xi}(\mathbf{t})=\mathbb{C}ov(\xi(\mathbf{t}),\xi(\mathbf{0})).$$
Assume that the covariance function satisfies the
following conditions:

\textbf{Assumption E1:} There exists a non-degenerate matrix $\mathbb{C}$ such that
$$r_{\xi}(\mathbb{C}\mathbf{t})=1-|t_{1}|^{\alpha_{1}}-|t_{2}|^{\alpha_{2}}+o(|t_{1}|^{\alpha_{1}}+|t_{2}|^{\alpha_{2}})$$
as $\mathbf{t}\to 0$ with $\alpha_{i}\in(0,2]$;

\textbf{Assumption E2:}
$r_{\xi}(\mathbf{t})<1$ for $\mathbf{t}\neq \mathbf{0}$.

To state two key lemmas, we recall the following type of Pickands constant.
For constant $a>0$, let
$$\mathcal{H}_{\alpha}^{a}[0,\lambda]=E\exp\left(\max_{ak\in[0,a\lambda]}\sqrt{2}B_{\alpha}(ak)-(ak)^{\alpha}\right)$$
and define
$$\mathcal{H}_{\alpha}(a)=\lim_{\lambda\rightarrow\infty}\frac{\mathcal{H}_{\alpha}^{a}[0,\lambda]}{a\lambda}.$$

We need the following results for the proofs of our main results.

\BL\label{Le:A61} Let $q_{i}=du^{-2/\alpha_{i}}$ for some $d>0$ and assume that  $\textbf{E1}$ and  $\textbf{E2}$ hold.
Then for any fixed rectangle $\mathbf{I_{h}}=[0,h_{1}]\times[0,h_{2}]$, we have
$$P\left(\max_{\mathbf{t}\in \mathbf{I_{h}}}\xi(\mathbf{t})> u\right)=
 h_{1}h_{2}\mathcal{H}_{\alpha_{1}} \mathcal{H}_{\alpha_{2}}|det \mathbb{C}^{-1}|u^{2/\alpha_{1}+2/\alpha_{2}}\Psi(u)(1+o(1))$$
and
$$P\left(\max_{\mathbf{kq}\in \mathbf{I_{h}}}\xi(\mathbf{kq})> u\right)=
 h_{1}h_{2}\mathcal{H}_{\alpha_{1}}(d) \mathcal{H}_{\alpha_{2}}(d)|det \mathbb{C}^{-1}|u^{2/\alpha_{1}+2/\alpha_{2}}\Psi(u)(1+o(1))$$
 as $u\rightarrow\infty$. The results also hold for the case $h_{1}=u^{-2/\alpha'}$ and $h_{2}=u^{-2/\alpha''}$ for $\alpha'>\alpha_{1}$ and $\alpha''>\alpha_{2}$.
\EL

\textbf{Proof:}  The first and second assertions can be proved following the proof of Lemma 7.1 of Piterbarg (1996) with some obvious changes, see also
the proof of Lemma 1 of D\c{e}bicki,  Hashorva and Soja-Kukiela (2015).
The third assertion follows from the proofs of the former two by using the double sums method, see the proof of Theorem 7.2 of Piterbarg (1996). $\Box$

\BL\label{Le:A62} Let $q_{i}=du^{-2/\alpha_{i}}$ for some $d>0$ and choose two constants $\alpha'>\alpha_{1}$ and $\alpha''>\alpha_{2}$.
Assume that  $\textbf{E1}$ and  $\textbf{E2}$ hold.
Then for the rectangle $\mathbf{I}=[0,u^{-2/\alpha'}]\times[0,u^{-2/\alpha''}]$, we have
$$P\left(\max_{\mathbf{kq}\in \mathbf{I}}\xi(\mathbf{kq})\leq u\right)
-P\left(\max_{\mathbf{t}\in \mathbf{I}}\xi(\mathbf{t})\leq u\right)
\leq |det \mathbb{C}^{-1}|\rho(d)u^{2/\alpha_{1}+2/\alpha_{2}-2/\alpha'-2/\alpha''}\Psi(u),$$
where $\rho(d)=\mathcal{H}_{\alpha_{1}}(d)\mathcal{H}_{\alpha_{2}}(d)-\mathcal{H}_{\alpha_{1}}\mathcal{H}_{\alpha_{2}}\rightarrow0$ as $d\rightarrow0$.
\EL

\textbf{Proof:} It is an immediate consequence of Lemma 6.1. $\Box$

\bigskip
{\bf Acknowledgement}: The author would like to thank Professor Enkelejd Hashorva for several suggestions and discussions.
The author also would like to deeply thank  the referees and the Associate Editor for useful comments and corrections which improved this paper
significantly.

\end{document}